\documentclass[12pt]{amsart}
\usepackage{amsmath, amsthm, amssymb}
\usepackage{fullpage}
\usepackage{xcolor}
\usepackage{hyperref}
\usepackage{soul}
\usepackage{enumitem}
\usepackage[makeroom]{cancel}

\newtheorem{theorem}{Theorem}
\newtheorem{lemma}{Lemma}

\newtheorem{corollary}[lemma]{Corollary}
\newtheorem{definition}[lemma]{Definition}

\newtheorem{conjecture}{Conjecture}
\numberwithin{lemma}{section}

\numberwithin{equation}{section}

\newcommand{\R}{{\mathbb R}}
\newcommand{\C}{{\mathbb C}}
\newcommand{\N}{{\mathbb N}}

\renewcommand{\R}{\mathbb R}

\newcommand{\bM}{\mathbf M}
\newcommand{\bP}{\mathbf P}
\newcommand{\bE}{\mathbf E}

\newcommand{\bu}{{\bar u}}

\newcommand{\la}{\langle}
\newcommand{\ra}{\rangle}

\newcommand{\calR}{\mathcal{R}}

\begin{document}

\title{The global well-posedness conjecture for  1D cubic dispersive equations}

\author{Mihaela Ifrim}
\address{Department of Mathematics, University of Wisconsin, Madison}
\email{ifrim@wisc.edu}

\author{ Daniel Tataru}
\address{Department of Mathematics, University of California at Berkeley}
\email{tataru@math.berkeley.edu}

\begin{abstract}

The goal of this article is to discuss 
a recent conjecture of the two authors,
which aims to describe the long time 
behavior of solutions to one-dimensional dispersive equations with cubic and higher nonlinearities. These problems  arguably represent  the single most important example where, 
even for small initial data,
the nonlinear effects are stronger than the dispersive effects. Consequently, 
the outcome predicted by the conjecture
depends essentially on the structure of the nonlinearity, precisely its focusing or defocusing character.
\end{abstract}


\maketitle

\setcounter{tocdepth}{1}
\tableofcontents


\section{Introduction}

\subsection{Linear vs nonlinear effects in dispersive flows}

A key comparison in the study of long time dynamics for nonlinear dispersive equations is between the relative strength of linear  and nonlinear effects. One naturally distinguishes three scenarios, which we briefly discuss below in the context of one-dimensional flows. As a good example 
we consider the nonlinear Sch\"odinger  problem  (NLS) in $\R \times \R$
\begin{equation}\label{nls-p}
iu_t + \Delta u = \pm u |u|^{p-1}, \qquad u(0) = u_0,   
\end{equation}
where typically the initial data is taken in Sobolev spaces $H^s(\R)$. 

\medskip

\textbf{(i) linear effects are stronger.} Then the expected behavior
for small initial data is \emph{scattering}, which means
that as $t \to \infty$ the solutions approach a solution for the associated linear flow. In the context of \eqref{nls-p},
in one space dimension this corresponds to $p > 5$, and global well-posedness and scattering may be easily proved using Strichartz estimates.

\medskip

\textbf{(ii) linear and nonlinear effects are balanced.} Then the expected behavior for small initial data is the same as in case (i). In the context of \eqref{nls-p}
in one space dimension this corresponds to $p = 5$, and global well-posedness and scattering may be again easily proved using Strichartz estimates.

\medskip

\textbf{(iii) nonlinear effects are stronger.} Here scattering 
cannot hold in general for $H^s$ data, and the question of global in time behavior of small data solutions is largely open. In the context of \eqref{nls-p}, in one space dimension this corresponds to $p < 5$. An interesting sub-problem here is the case of small and localized data, where the prospects of scattering do somewhat improve.

\medskip

In this note we consider the single most important example 
where the nonlinear effects are stronger,
namely the case of one-dimensional dispersive flows with cubic nonlinearities.
This corresponds to $p = 3$ in \eqref{nls-p}. The cubic power is interesting here for multiple reasons:
\begin{itemize}

\item $3$ is the only odd integer in the corresponding range of exponents, i.e. below $5$.

\item cubic problems represent universal models in nonlinear 
dispersive equations.

\item $p=3$ represents the threshold for scattering for small and localized data. Precisely, in this context scattering holds for $p > 3$, but only a modified form of scattering occurs 
at $p = 3$, as explained in Section~\ref{s:localized}. 
\end{itemize}

\subsection{ Cubic dispersive flows}
For later reference, the reader should keep in mind a model of the form 
\begin{equation}\label{cubic}
i u_t - A(D) u = C(u,\bar u,u)   \qquad \text{ in } 
\R \times \R,
\end{equation}
with initial data 
\begin{equation}
u(0) = u_0 \in H^s(\R)  ,
\end{equation}
where we make some reasonable assumptions:
\begin{itemize}
    \item The dispersion relation $\xi \to a(\xi)$ 
is given by a smooth, real symbol $a$ with nondegenerate  
dispersion, i.e. $a''(\xi)  \neq 0$.

\item The cubic nonlinearity $C (\cdot\, ,\,  \cdot \, ,\, \cdot )$ is translation invariant, and thus determined by its symbol 
$c(\xi_1,\xi_2,\xi_3)$.  The choice of variables $(u,\bar u,u)$ guarantees phase rotation symmetry, which is a standard property for many flows of this type.

\item Quintic and higher terms may also be added, but do not affect the discussion below.
\end{itemize}

From a local well-posedness perspective  one may distinguish two scenarios:

\begin{enumerate}[label=\roman*)]
\item semilinear, where solutions have a Lipschitz continuous dependence on the data.
\item quasilinear, where solutions only have a continuous dependence on the data.
\end{enumerate}
These two scenarios continue to differ significantly when examined from a global well-posedness standpoint. 

\subsection{ The long time behavior conjectures}

While a large array of global well-posedness results have been proved for cubic 1D problems in the case of small and localized data, until last year little was known concerning global dynamics for data which is merely small in Sobolev spaces $H^s$.  This is when the two authors formulated 
a broad conjecture which applies to this 
problem. Our global well-posedness (GWP) conjecture, which applies to both semilinear and quasilinear 1D problems, differentiates between focusing and defocusing  problems, as follows:

\begin{conjecture}[Non-localized data defocusing GWP conjecture \cite{IT-global}]\label{c:defocusing}
One-dimensional dispersive flows on the real line with cubic defocusing nonlinearities and small initial data have global in time, ``dispersive" solutions.
\end{conjecture}

Critically, compared with any earlier work
in this direction, this conjecture requires 
no localization for the initial data.
``Dispersive"  here is interpreted in a weak sense, to mean that the solution satisfies global $L^6_{t,x}$ Strichartz estimates and bilinear $L^2_{t,x}$ bounds. This is due to the strong nonlinear effects, which preclude classical scattering for any solutions to such problems. This is discussed in greater detail in Section~\ref{s:conjecture}.

The defocusing condition for the nonlinearity is essential in the GWP conjecture. In the focusing case, the existence of small amplitude solitons generally prevents global, scattering solutions. Nevertheless, in another recent paper,  the authors
have conjectured that instead, in the focusing case long time solutions can be  obtained on a likely optimal time-scale: 

\begin{conjecture}[Non-localized data long time well-posedness conjecture \cite{IT-focusing}]\label{c:focusing}
One-dimensi\-o\-nal \mbox{dispersive} flows on the real line with cubic conservative nonlinearities and  initial data of size $\epsilon \ll 1$ have  long time solutions on the $\epsilon^{-8}$ time-scale.
\end{conjecture}

Here there is no need to explicitly assume that the problem is 
focusing. Instead we impose a weaker, \emph{conservative} assumption, which is very natural, and heuristically aims to prevent nonlinear ode blow-up of solutions with wave-packet localization; see the discussion in Section~\ref{s:wp}.  This assumption is implicitly satisfied 
for the defocusing problems in Conjecture~\ref{c:defocusing}.
The solutions which are the subject of the second conjecture will also be expected
to satisfy $L^6_{t,x}$ Strichartz estimates and bilinear $L^2_{t,x}$ bounds on suitable time scales;  this 
is  discussed  in Section~\ref{s:conjecture}.

Our aim in this paper is to provide context and motivation for these conjectures, leading to a more precise formulation,
as well as a summary  of our progress up to this point, and an overview of  the novel approach we have developed in order to prove these results.

\subsection{ An overview} This paper is organized as follows:

\medskip

\textbf{1. Linear dispersive estimates.} Estimates for solutions to a linear 1D dispersive evolution
are discussed in Section~\ref{s:linear}. This includes 
dispersive decay, Strichartz bounds as well as bilinear $L^2_{t,x}$ bounds. While these estimates are never used for the nonlinear problem due to the nonperturbative nature of the nonlinearity on large time-scale, they still serve as a guide for what can be expected.

\medskip

\textbf{2. Wave packet dynamics.} One critical 
contribution to our intuition regarding these conjectures
is gained by considering wave packet dynamics in Section~\ref{s:wp}.
This shows that  for wide wave packets  the cubic nonlinear effects are seen before the linear ones, which in turn  indicates that no scattering can be expected
for such problems. The role played by solitary waves 
in the focusing case is also described there.

\medskip

\textbf{3. The special case of small and localized data} 
was previously considered in many earlier works. 
As a prelude to our main discussion, this case is  
briefly discussed in Section~\ref{s:localized}. One reason
this case is interesting is that it is the only case 
where some form of \emph{modified scattering} may be expected. 

\medskip

\textbf{4. The case of small but nonlocalized data}
is the setting of our conjectures. These are discussed in greater detail in Section~\ref{s:conjecture}. In particular, we provide there a more precise form of the conjectures,
which is based both on the earlier heuristic discussion 
and on the positive results so far.

\medskip

\textbf{5. Validating the conjectures.} The goal of our work so far has been to prove the above conjecture 
for certain well-chosen models. Our results are presented in Section~\ref{s:results}, and cover two settings: 
\begin{itemize}
\item Semilinear Schr\"odinger flows, which were considered in  \cite{IT-global}, \cite{IT-focusing}. 
\item Quasilinear Schr\"odinger flows, which were considered in  \cite{IT-qnls}.
\end{itemize}

\medskip

\textbf{6. The key ideas of our approach.}
These are presented in the last section of the 
article, in the context of the results in \cite{IT-global}, \cite{IT-focusing}, \cite{IT-qnls}. But we expect 
the same ideas to work in a much broader context.

\subsection{Acknowledgements}
The first author was supported by the Sloan Foundation, and by an NSF CAREER grant DMS-1845037. The second author was supported by the NSF grant DMS-2054975 as well as by a Simons Investigator grant from the Simons Foundation.

\section{Linear dispersive decay} \label{s:linear}

In this section we discuss the dispersive 
decay properties for the corresponding linear flow
\begin{equation}\label{linear}
i u_t -  A(D) u = 0,  \qquad
u(0) = u_0  .
\end{equation}

\subsection{The fundamental solution 
and dispersive estimates} Here we begin with the fundamental solution, and then consider the  corresponding uniform decay 
properties for solutions with localized data.
The fundamental solution is given by the group  of $L^2$ isometries 
\[
S(t) = e^{- i t A(D)},
\]
and its kernel is given by the inverse Fourier transform of the symbol,
\[
K(t,x) = \frac{1}{2\pi} \int_\R e^{-ita(\xi)+x\xi} \, d\xi,
\]
which is interpreted as an oscillatory integral. This can be evaluated using the method of stationary phase, which,
assuming say that $a$ is strictly convex and coercive, yields an asymptotic expansion for the fundamental solution of the form
 \begin{equation}
    K(t,x) \approx  \frac{1}{\sqrt{2 \pi t a''(\xi_v)} } e^{-\frac{i\pi}{4}} e^{i t \phi(v)} + O(t^{-1}), \qquad v = x/t,
  \end{equation}  
 where the phase function $\phi$ is the Legendre 
 transform of $a$, precisely 
\[
\phi(v) = \sup_{\xi \in \R} \left\{ v \xi - a(\xi)\right\},
\]
and $v$ and $\xi_v$  satisfy the relations 
 \[
 a'(\xi_v) = v, \qquad \phi'(v) = \xi_v.
 \]
Here, $v$ represents the group velocity of waves with frequency $\xi_v$.
 
It is interesting to consider the behavior of $K$ 
along rays $x = vt$, which can be thought of as  the trajectories of propagating wave originating from (near) $x = 0$.  We immediately see that 
the fundamental solution has $t^{-\frac12}$ decay, 
which we will refer to as the dispersive decay on one space dimension.

A corollary of this decay property for the fundamental solution is a similar decay property for general solutions
for the linear equation \eqref{linear} with ``nice'' and localized (e.g. Schwartz) initial data. Even better, for such solutions 
we obtain an asymptotic expansion
\begin{equation}\label{lin-asympt}
u(t,x) =    \frac{\gamma(v)}{\sqrt{ta''(\xi_v)}} e^{i t \phi(v)} + O(t^{-1}).
\end{equation}
The function $\gamma$  can  be thought as an \emph{asymptotic profile} for the solution $u$,
which in the linear case is directly related to the 
Fourier transform of the initial data,
\[
\gamma(v) =  e^{-\frac{i\pi}4} \hat u_0(\xi_v).
\]
One may ask whether  such uniform decay properties may hold for the corresponding nonlinear problem \eqref{cubic}. At best,
this would also require the initial data to be localized. 
However, we will see in Section~\ref{s:localized} that 
even then this cannot happen in the case of cubic nonlinearities.

\subsection{Strichartz and bilinear $L^2_{t,x}$ bounds} \label{s:disp}

Both of these bounds will be the type of estimates 
we expect later on for the solutions to both 
the full nonlinear problem and for its linearization. For reference and comparison purposes, here we briefly recall these bounds in the case of the constant coefficient
Schr\"odinger flow.

We begin with the classical  Strichartz inequality, which applies to solutions to the inhomogeneous linear Schr\"odinger equation:
\begin{equation}\label{lin-inhom}
(i\partial_t + \partial^2_x)u = f, \qquad u(0) = u_0.
\end{equation}

To measure the solution $u$  we will use the Strichartz space $S$ associated to the $L^2$ flow,  defined by 
\[
S := L^\infty_t L^2_x \cap L^4_t L^\infty_x.
\]
For the source term $f$ we will use the  dual Strichartz space
\[
S' = L^1_t L^2_x + L^{\frac43} _t L^1_x .
\]

The Strichartz estimates in the $L^2$ setting are a consequence of the nonvanishing curvature 
of the characteristic set for the Schr\"odinger equation, i.e. the parabola $\{\tau +\xi^2 = 0\}$.
They are summarized in the following:

\begin{lemma}\label{l:Strichartz}
Assume that $u$ solves \eqref{lin-inhom} in $[0,T] \times \R$. Then 
the following estimate holds.
\begin{equation}
\label{strichartz}
\| u\|_S \lesssim \|u_0 \|_{L^2} + \|f\|_{S'} .
\end{equation}
\end{lemma}
One intermediate norm between the two endpoints in $S$ is $L^6_{t,x}$, and its dual is $L^\frac65 _{t,x}$. In our estimates later in the paper, we will give preference to this Strichartz norm and neglect the rest of the family.

\bigskip

The second property of the linear Schr\"odinger equation we want to describe here is the 
bilinear $L^2_{t,x}$ estimate. This applies for the product of two  waves with separated Fourier supports, and thus with separated group velocities. Because of this, such bilinear bounds are better interpreted as transversality estimates, rather than Strichartz estimates. The precise statement  is as follows:

\begin{lemma}
\label{l:bi}
Let $u^1$, $u^2$ be two solutions to the inhomogeneous Schr\"odinger equation \eqref{lin-inhom} with data 
$u^1_0$, $u^2_0$ and inhomogeneous terms $f^1$ and $f^2$. Assume 
that the sets 
\[
E_i =  \text{supp } \hat u^i
\]
are disjoint. Then we have 
\begin{equation}
\label{bi-di}
\| u^1 u^2\|_{L^2} \lesssim \frac{1}{\text{dist}(E_1,E_2)^\frac12} 
( \|u_0^1 \|_{L^2} + \|f^1\|_{S'}) ( \|u_0^2 \|_{L^2} + \|f^2\|_{S'}).
\end{equation}
\end{lemma}

One corollary of this applies in the case when we look at the product
of two solutions which are supported in different dyadic regions:
\begin{corollary}\label{c:bi-jk}
Assume that $u^1$ and $u^2$ as above are supported in dyadic regions $\vert \xi\vert \approx 2^j$ and $\vert \xi\vert \approx 2^k$, $\vert j-k\vert >2$, then
\begin{equation}
\label{bi}
\| u^1 u^2\|_{L^2} \lesssim 2^{-\frac{\max\left\lbrace j,k \right\rbrace  }{2}}
( \|u_0^1 \|_{L^2} + \|f^1\|_{S'}) ( \|u_0^2 \|_{L^2} + \|f^2\|_{S'}).
\end{equation}
\end{corollary}

Another useful  case is when we look at the product
of two solutions which are supported in the same  dyadic region, but with frequency separation:
\begin{corollary}\label{c:bi-kk}
Assume that $u^1$ and $u^2$ as above are supported in the dyadic region 
$\vert \xi\vert \approx 2^k$, but have $O(2^k)$ frequency separation between their supports.
Then 
\begin{equation}
\label{bi-kk}
\| u^1 u^2\|_{L^2} \lesssim 2^{-\frac{k}{2}}
( \|u_0^1 \|_{L^2} + \|f^1\|_{S'}) ( \|u_0^2 \|_{L^2} + \|f^2\|_{S'}).
\end{equation}
\end{corollary}

\section{Heuristic considerations via wave packet dynamics} \label{s:wp}

\subsection{ Resonant analysis} In this subsection we introduce the notion of \emph{resonant frequencies}, which applies to nonlinear interactions relative to the linear dispersion relation associated to our equation. Given our setup, we limit the discussion to the case of cubic nonlinearities. We also assume that our problem has the phase rotation symmetry, which implies that the cubic part $C$ of the nonlinearity can be thought of as a trilinear form with arguments $C(u,\bu,u)$.

For such a trilinear form, given three input frequencies $\xi_1, \xi_2,\xi_3$ for 
our cubic nonlinearity, the output will be at frequency 
\[
\xi_4 = \xi_1-\xi_2+\xi_3.
\]
This relation can be described in a more symmetric fashion as 
\[
\Delta^4 \xi = 0, \qquad \Delta^4 \xi := \xi_1-\xi_2+\xi_3-\xi_4 .
\]
This is a resonant interaction if and only if we have a similar relation for the associated time frequencies, namely 
\[
\Delta^4 \xi^2 = 0, \qquad \Delta^4 \xi^2 := \xi_1^2-\xi_2^2+\xi_3^2-\xi_4^2 .
\]
Hence, we define the resonant set in a symmetric fashion as 
\[
\calR := \{ \Delta^4 \xi = 0, \ \Delta^4 \xi^2 = 0\}.
\]
It is easily seen that this set may be characterized as 
\[
\calR = \{ \{\xi_1,\xi_3\} = \{\xi_2,\xi_4\}\}.
\]

When considering estimates for cubic resonant interactions, it is clear that the case when $\xi_1 \neq \xi_3$ is more favourable, as there 
we have access to bilinear $L^2_{t,x}$ bounds. 
The unfavourable case is when all four frequencies are equal. We denote this set by 
\[
\calR_2 = \{\xi_1 =\xi_3 = \xi_2=\xi_4\},
\]
and we will refer to it as \emph{the doubly resonant set}. Heuristically, this set carries the bulk 
of the cubic long range interactions in the nonlinear problem, and can be associated with wave packet self-interactions.

\subsection{Linear wave packets}

Wave packets are the most concentrated solutions
to dispersive equations, and also those solutions for which approximately optimize the Strichartz estimates.
To streamline the following discussion, let us assume
that we are looking at solutions to the linear equation 
\eqref{linear} which are frequency localized in a compact set. Within this set, we consider a center frequency $\xi_0$
and a frequency scale $\delta \xi = N \ll 1$.
Correspondingly, we take the initial data $u_0$ to be localized in this set. Then the solution will retain this localization, and at leading order will move with group velocity $a'(\xi_0)$. Precisely, the transport in this direction with the appropriate phase yields an approximate solution
\[
u^{app}(t,x) = u_0 (x-ta'(\xi_0)) e^{it(\xi_0 a'(\xi_0)-a(\xi_0)) },
\]
where the error $f$ given by
\[
(i \partial_t - A(D))u^{app} = (a(\xi_0) -(D +\xi_0) a'(\xi_0) - A(D)) u_0 (x-ta'(\xi_0)) e^{it(\xi_0 a'(\xi_0)-a(\xi_0)) } :  = f
\]
satisfies
\[
\|f(t) \|_{L^2} \lesssim N^{2} \|u_0\|_{L^2}.
\]
This implies that $u^{app}$ can be seen as a good approximation to the linear flow on the time-scale
\begin{equation}
\delta t^{lin} = N^{-2}.    
\end{equation}

So far we have only considered the frequency localization of the data. However, in order to obtain spatial concentration for the solutions one should also work 
with spatial localization, which by the uncertainty principle can at best be achieved on the $N^{-1}$ scale.
A good setup here is to take the initial data $u_0$ of the form
\begin{equation}\label{wp-data}
u_0(x) = N^{\frac12} \phi_0(N(x-x_0)) e^{i x \xi_0}
\end{equation}
with $\phi_0$ a fixed Schwartz function,
and the $N^\frac12$ factor added for $L^2$ normalization.
We then represent the solution $u$ as 
\begin{equation}\label{wp-ansatz}
u(t,x) = N^{\frac12} \phi (t,N(x-x_0 -ta'(\xi_0))) e^{ix\xi_0}  e^{it(\xi_0 a'(\xi_0)-a(\xi_0))} .
\end{equation}
Following the above heuristic, the function $\phi$ 
will remain uniformly in the Schwartz space on the 
time interval $\{|t| \lesssim  \delta t^{lin} = N^{-2}\}$.
In effect the function $\phi$ may be shown to solve 
an equation which is similar to \eqref{linear} but with 
$A$ replaced by $N^2 \tilde A$ where 
\[
\tilde a(\xi) =  N^{-2}[a(\xi_0+N\xi) - a(\xi_0) - N \xi a'(\xi_0)\xi_0]\ ,
\]
and $\tilde a(0) = 0$ and $\tilde a'(0) = 0$.
This corresponds to composing a Galilean transformation with a spatial rescaling. 

We call such solutions \emph{wave packets}, noting that they are spatially localized in an interval of size $N^{-1}$ moving 
with velocity $a'(\xi_0)$ for a time $O(N^{-2})$.

\subsection{Trilinear wave packet self-interactions}
Here we turn our attention to the nonlinear problem \eqref{cubic}, where we start with initial data exactly 
as the wave packets above.  Then the function $\phi$ 
also solves an equation of the form \eqref{cubic}, precisely 
\begin{equation}
i \phi_t - N^2 \tilde A(D) \phi = N \tilde C(\phi,\bar \phi,\phi)   \qquad u(0) = u_0,
\end{equation}
where $\tilde A$ is as in the linear case, while
\[
\tilde c(\xi_1,\xi_2,\xi_3) =
c(\xi_0+N \xi_1,\xi_0+N \xi_2, \xi_0+N \xi_3).
\]
Here on the unit frequency scale we have 
\[
\tilde c(\xi_1,\xi_2,\xi_3) =
c(\xi_0,\xi_0, \xi_0)+O(N).
\]
Therefore we can write the above equation as 
\begin{equation}
i \phi_t - N^2 \tilde A(D) \phi = N c(\xi_0,\xi_0,\xi_0)
\phi |\phi|^2
+ O(N^2).
\end{equation}
The $N^2$ terms will only be effective on the $N^{-2}$ 
time-scale, but the $N$ term will already be visible 
on the $N^{-1}$ time-scale. Isolating this term we obtain the \emph{reduced equation}
\begin{equation}\label{reduced}
i \phi_t = N c(\xi_0,\xi_0,\xi_0)
\phi |\phi|^2,
\end{equation}
which is an ode. where the $N$ factor may be removed by a 
time rescaling $Nt = s$. For this ode we distinguish 
two scenarios:

\begin{enumerate}
    \item $\Im c(\xi_0,\xi_0,\xi_0) \neq  0$. Then 
finite time blow-up at the $N^{-1}$ time-scale may happen even for small initial data. 

\item $\Im c(\xi_0,\xi_0,\xi_0) = 0$. In this case the ode dynamics reduce to a phase rotation,
\[
\phi(t,x) = \phi_0(x) e^{-i N t  c(\xi_0,\xi_0,\xi_0) |\phi_0(x)|^2},
\]
where the solutions remain globally bounded. It is this case that is of interest to us from the perspective of obtaining long time or global solutions. Later we will refer to this case as \emph{conservative}, see Definition~\ref{d:conservative}.
\end{enumerate}

In the first case the ode blow-up in the reduced equation does not immediately  imply blow-up also for the full equation, but does indicate that unbounded norm growth will occur on the same time-scale. Precisely, if the 
initial data $\phi_0$ is assumed to have size $\epsilon \ll 1$, then the solution can be shown to reach size $1$ 
by the time $N^{-1} \epsilon^{-2}$.

In the second case the reduced model dynamics are global and bounded, but they are already effective on the $N^{-1}$ time-scale, which is below the $N^{-2}$ linear time-scale. Further, we note that beyond the $N^{-1}$ time-scale the phase rotation effects become substantial 
since $|\phi_0|^2$ is varying. Precisely, for $t > N^{-1}$ the solution $\phi$ is localized at frequencies $\lesssim Nt$, which is no longer consistent with the 
unit localization scale implicit in our ansatz. At this point it will be the focusing and defocusing effects
which take precedence; these in turn depend on the sign of $c(\xi_0,\xi_0,\xi_0)$ relative to the sign of $a''(\xi_0)$.

\subsection{Solitons} These are special solutions for \eqref{cubic} which have the form 
\begin{equation}
u(x) = \phi(x-x_0-vt) e^{i \omega t} 
\end{equation}
with a time independent profile $\phi$ traveling with velocity $v$ and a time dependent phase 
rotation with velocity $\omega$. Such solutions are reasonable to consider in the case when the problem has a phase rotation symmetry. In that case, by choosing the frequency $\xi_0$ so that the associated group velocity is $v = a'(\xi_0)$,
and the soliton scale is $N$, one may write an equation for 
\[
\tilde \phi (x) = N e^{-i(x-x_0) \xi_0} \phi(N x).
\]
This has the form 
\begin{equation}\label{soliton}
\tilde A (D) \tilde \phi = \tilde C(\tilde \phi, \bar{\tilde \phi}, \tilde \phi) - \tilde \omega \tilde \phi     , 
\end{equation}
where $\tilde\omega = N^{-2} \omega$. 
Assuming for normalization purposes that $c(\xi_0,\xi_0,\xi_0) = \pm 2$, this equation is well approximated as $N \to \infty$ by  the cubic NLS soliton equation
\begin{equation}\label{nls3-soliton}
- \Delta  Q = \pm 2 Q|Q|^2  - \omega Q.
\end{equation}
This is known to admit a localized solution only in the focusing case, i.e. the ``+" sign, and with $\omega > 0$. For $\omega = 4$ this has the form
\begin{equation}\label{exact-soliton}
Q = 2 \mathop{sech} 2x,
\end{equation}
and for other values it is obtained by rescaling.

Under reasonable assumptions, a perturbation argument can be used to show that in the focusing case \eqref{soliton}
has a solution $\tilde \phi$ close to above the $NLS^3$ soliton. The corresponding soliton profile $\phi$ 
in \eqref{soliton} has the following properties:

\begin{itemize}
    \item Localization: it is localized around position $x_0$ on the $N^{-1}$ scale and around frequency $\xi_0$ on the $N$ scale.
    \item Size: it has $L^2$ size $O(N^\frac12)$.
\end{itemize}
As $N$ is an arbitrarily small parameter, we will refer 
to such solitons as \emph{small solitons}, i.e. which can arise from arbitrarily small $H^s$ data.

This implies that in general we cannot expect 
to have global dispersive solutions in the focusing 
case. Pursuing this further, the time-scale in our focusing conjecture is motivated precisely by these 
solitons, which maximize the Strichartz norms in our 
bounds for the solutions.

\subsection*{Conclusions} We summarize the conclusions of our heuristic considerations above as follows:

\begin{enumerate}[label=(\roman*)]
\item  The condition $c(\xi,\xi,\xi) \in \R$ is necessary 
in order to prevent ode growth/blow-up phenomena.

\item  At low frequency scales nonlinear wave packet solutions are no longer well approximated by linear wave 
packet solutions, so no linear scattering is possible in 1D cubic dispersive flows.

\item  Focusing/defocusing effects are important in the long time behavior of the solutions.
\end{enumerate}

\section{The case of localized data and modified scattering} \label{s:localized}

For comparison purposes, here we discuss the much easier case of small, smooth and data, following our expository paper \cite{IT-packet}. A reasonable assumption would be say
\begin{equation}\label{small-loc}
\|u_0\|_{H^s_x}+ \| x u_0\|_{L^2_x} \leq \epsilon \ll 1.  \end{equation}
For such data, we ask whether the solution is global, and also we ask 
about its asymptotic behavior at infinity.

To simplify the discussion, it is convenient to assume that $a$ is also coercive,
so that the map $\xi \to v_\xi = a'(\xi)$ is surjective; but see \cite{IT-packet} for other cases.
We also set aside the question of minimizing the regularity where global well-posedness results hold, referring the reader again to \cite{IT-packet}.

The naive expected behavior for the global solutions solutions would be the one given by the linear expression \eqref{lin-asympt}, 
\begin{equation}\label{non-lin-try}
u(t,x) \approx    \frac{\gamma(v)}{\sqrt{ta''(\xi_v)}} e^{i t \phi(v)} .
\end{equation}

However, assuming an ansatz 
of this type with a nice and localized profile $\gamma$ one may 
show that the linear part of the equation 
\eqref{cubic} has size
\[
(i \partial_t - A(D) )u = O(t^{-2}), 
\]
whereas the nonlinear part is
\[
C(u,\bar u, u) = \frac{c(\xi_v,\xi_v,\xi_v) \gamma(v)|\gamma(v)|^2}{(t a''(\xi_v))^{\frac32}} e^{it \phi(v)} + O(t^{-2})
\]
The two expressions are mismatched, so our ansatz
\eqref{non-lin-try} cannot be valid.  However, as it turns out, we are not very far from reality. 
Heuristically, this is because the initial data's spatial localization both prevents localization on frequency scales $\delta \xi \ll 1$, and forces spatial separation between the long time evolution of different frequencies. One may also calculate 
the $L^2$ norm of the error 
\[
e(t) = (i \partial_t - A(D) )u - C(u,\bar u, u)
\]
at this stage. This gives
\[
\|e(t)\|_{L^2} \approx t^{-1}.
\]
This fails to be integrable in time, but 
it is integrable on dyadic time-scales. 
One may interpret this fact as indicating that our 
initial ansatz is accurate on dyadic time-scales, even if not globally in time. To rectify this, it is natural to slowly modulate our asymptotic profile 
$\gamma$, allowing it to vary on a dyadic time-scale. Therefore, we substitute the ansatz \eqref{non-lin-try} with a more flexible alternative:
\begin{equation}\label{non-lin}
u(t,x) \approx    \frac{\gamma(t,v)}{\sqrt{ta''(\xi_v)}} e^{i t \phi(v)} .
\end{equation}
This leaves the approximate expression for $C(u,\bar u, u)$ unchanged, but the linear part acquires a new leading term
\[
(i \partial_t - A(D) )u =  \frac{i\partial_t \gamma(t,v)}{\sqrt{ta''(\xi_v)}} e^{i t \phi(v)} + O(t^{-2}). 
\]
Matching this with the corresponding linear part
we arrive at an \emph{asymptotic equation}
for the \emph{asymptotic profile} $\gamma$ along rays, namely
\begin{equation}
 i\partial_t \gamma(t,v) =     \frac{c(\xi_v,\xi_v,\xi_v)}{t a''(\xi_v)}  \gamma(t,v)|\gamma(t,v)|^2.
\end{equation}
This is an ode whose solutions remain of constant size, provided that the diagonal trace $c(\xi,\xi,\xi)$ is real valued. Further, the solutions can be represented as
\[
\gamma(t,v) = \gamma_0(v) e^{i \log t \ b(v) |\gamma_0(v)|^2}, \qquad b(v) =  \frac{c(\xi_v,\xi_v,\xi_v)}{a''(\xi_v)}
\]
for a well chosen radiation profile $\gamma_0$.
For convenience we state a typical theorem, which holds under suitable assumptions on $a$ and $c$.

\begin{theorem}
Assume that the initial data $u_0$ for \eqref{non-lin} satisfies \eqref{small-loc}. Then the solution $u$ is global in time, and admits an asymptotic expansion of the form
\begin{equation}\label{non-lin+}
u(t,x) =   \frac{\gamma_0(v)}{\sqrt{ta''(\xi_v)}}  e^{i \log t \ b(v) |\gamma_0(v)|^2} e^{i t \phi(v)} + O(t^{-1-}) .
\end{equation}
\end{theorem}
Results of this type have been proved for many cubic dispersive equations in one space dimension, both 
in the semilinear and in the quasilinear setting.
For more details, we refer the reader to \cite{IT-packet}, but for a few additional reference samples see also \cite{HN,HN1,LS,KP,IT-NLS}. Asymptotic expansions of the form \eqref{non-lin+} are referred to as \emph{modified scattering}. 

\subsection*{Conclusions}
We complete our discussion of this case
with several remarks connecting this simpler setting with our setup for these notes:

\begin{enumerate}[label=(\roman*)]
    \item Localizing the initial data reduces the strength of nonlinear interactions to the point where it is only borderline nonperturbative. But we still get modified scattering rather than regular scattering. 

    \item The conservative assumption $c(\xi,\xi,\xi) \in \R$ is equally needed in the localized data case.

    \item The distinction between the focusing and defocusing scenarios is no longer needed here.
    This is because the small solitons are incompatible with the small and localized initial data.
    
\end{enumerate}

\section{Nonlocalized data: a more precise form of the conjectures}  \label{s:conjecture}

We write a general one dimensional dispersive flow
\begin{equation}\label{full}
i u_t - A(D) u = N(u,\bar u),  \qquad
u(0) = u_0  , 
\end{equation}
where in $N$ we place all the terms which are at least quadratic in $u$. Implicitly we assume here
that the nonlinearity $N$ is smooth as a function 
of $u$ and thus admits a Taylor series expansion around $u = 0$. 

From a linear perspective, we call such a problem (strictly) dispersive if the symbol $a$, commonly called the \emph{dispersion relation}, is real and $a''(\xi) \neq 0$. To fix the signs later on 
we will assume for simplicity that a is strictly convex,
\begin{equation}\label{a-sign}
  a''(\xi) > 0, \qquad \xi \in \R.  
\end{equation}
Of course in specific problems the behavior of $a''$ at $\pm \infty$ also plays a role.
\medskip

Depending on the strength of $N$, one may classify such problems as 
\begin{enumerate}[label = (\roman*)]
    \item  semilinear, i.e. with a perturbative $N$ on small time-scales, and thus with a Lipschitz dependence of the solution on the initial data locally in time, or
    \item quasilinear, with an $N$ that belongs to the principal part of the equation, and merely continuous dependence of the solution on the initial data. 
\end{enumerate}
Our conjecture does not discriminate between the two cases, and in any case the long time dynamics are expected to be nonperturbative even in the semilinear case. 

A prerequisite for any global well-posedness result 
is to have local well-posedness. This could be very problem specific, particularly in the quasilinear case. For this reason, we do not make local well-posedness part of our discussion, and instead implicitly assume that the problems which are within the scope of our conjectures are locally well-posed. 

In order to provide a more precise form of the conjectures it is convenient to make several simplifying assumptions. We emphasize these are not all
required, but instead provide a simple baseline,
which may be later developed in multiple directions.
We now carefully review these assumptions. We begin with 

\begin{definition}
 We say that the equation \eqref{full} is \emph{cubic} if  $N$ is at least cubic at zero. 
\end{definition}

This excludes quadratic interactions, though we 
note that there are also many situations where 
quadratic interactions may be replaced by cubic ones 
via a normal form analysis.

The second assumption we make concerns the symmetries
of the equation \eqref{full}. As written it is already invariant with respect to translations.
In addition, we will require a second symmetry, which is very common for instance in the NLS realm:

\begin{definition}
We say that the equation \eqref{full}
has  \emph{phase rotation} symmetry if it is 
invariant with respect to the transformation 
$u \to u e^{i\theta}$ for $\theta \in \R$.
\end{definition}

One should view this as a simplifying   assumption rather than a fundamental one. Its primarily role is to  streamline the setup
and to allow for a clean formulation  of the conjectures which is free from technical assumptions.  For instance, it guarantees that only terms with odd homogeneity are present in $N$, beginning with cubic ones. At the cubic level, it allows only for trilinear forms $C(u,\bar u,u)$.
Precisely, we can expand $N$ as  
\[
N(u,\bar u) = C(u,\bar u,u)+N^{\geq 5} (u,\bar u),
\]
where the remainder $N^{\geq 5} (u,\bar u)$ 
contains only quintic and higher order terms.
The remaining assumptions only refer to the cubic term $C$, which can be described by its symbol 
$c(\xi_1,\xi_2,\xi_3)$ as 
\[
\widehat{C(u,\bar u,u)}(\xi) = \frac{1}{2\pi} 
\int_{\xi_1-\xi_2+\xi_3 = \xi}  \hat u(\xi_1) \bar {\hat u}(\xi_2) \hat u(\xi_3) \, d \xi_1 d \xi_2.
\]
Motivated by our wave packet heuristics,
the first condition on $C$ is as follows:

\begin{definition}\label{d:conservative}
We say that the equation \eqref{full}
is \emph{conservative} if  the cubic component of the nonlinearity satisfies 
\[
c(\xi,\xi,\xi), \, \partial_{\xi_j}
c(\xi,\xi,\xi) \in \R .
\]
\end{definition}
Compared with the earlier wave packet heuristics here we have added a second condition on the gradient of the symbol $c$ on the diagonal. This is primarily motivated by our positive results; however, it remains unclear to what degree it is actually necessary.

The last goal is to distinguish between the focusing and defocusing scenarios. 

\begin{definition}\label{d:defocusing}
We say that the equation \eqref{full}
is \emph{defocusing} if   $c(\xi,\xi,\xi)$ 
    is positive definite,
    \begin{equation}\label{defocusing}
    c(\xi,\xi,\xi) > 0, \qquad \xi \in \R.
    \end{equation}
If $c$ has the opposite sign then we call the problem focusing.    
\end{definition} 
We remark that the sign in \eqref{defocusing}
depends on the sign convention in \eqref{a-sign}.
In specific examples, it is plausible that the above 
sign condition may need to be supplemented by 
a quantitative bound from below depending on the behavior of $a$ at $\pm \infty$.

We are now ready to provide more precise forms 
of our conjectures under these more restrictive assumptions. Our defocusing conjecture is as follows:

\begin{conjecture}[Non-localized data defocusing GWP conjecture \cite{IT-global}]\label{c:defocusing+}
Consider a one dimensional dispersive flow as in \eqref{full} which is locally well-posed in some Sobolev space $H^s$. Assume that our equation has 
phase rotation symmetry, it is conservative and defocusing.  Then for any initial data $u_0$ 
so that 
\begin{equation}
\label{small-data}
\|u_0\|_{H^s_x} \leq \epsilon \ll 1,
\end{equation}
the solution $u$ is global in time, uniformly bounded,
\begin{equation}
\label{small-sln}
\|u\|_{L^\infty_t H^s_x} \lesssim  \epsilon,
\end{equation}
and dispersive, in the sense that it satisfies 
bounds similar to the corresponding linear flow as follows:
\begin{itemize}
\item global $L^6_{t,x}$ Strichartz bounds, and 
\item global bilinear $L^2_{t,x}$ bounds.
\end{itemize}
\end{conjecture}

Here the $L^6_{t,x}$ Strichartz bounds are akin to
the ones provided by Lemma~\ref{l:Strichartz},
with an appropriate number of derivatives depending on both the convexity of $a$, the positivity of $c$  
on the diagonal and the Sobolev exponent $s$.
Similarly, the bilinear $L^2_{t,x}$ bound is akin to 
Lemma~\ref{l:bi} with an appropriate number
of derivatives, and also with a bound from below 
on the size of the sets $E_j$ and their separation.

\bigskip

We conclude with the focusing conjecture, which is as follows:

\begin{conjecture}[Non-localized data focusing long time well-posedness conjecture]\label{c:focusing+}
Consider a one dimensional dispersive flow as in \eqref{full} which is locally well-posed in some Sobolev space $H^s$. Assume that our equation has 
phase rotation symmetry and is conservative.  Then for any initial data $u_0$ 
so that 
\begin{equation}
\label{small-data-focus}
\|u_0\|_{H^s_x} \leq \epsilon \ll 1,
\end{equation}
the solution $u$ exists for 
\[
|t| \ll \epsilon^{-8}
\]
is uniformly bounded,
\begin{equation}
\label{small-sln-focus}
\|u\|_{L^\infty_t H^s_x} \lesssim  \epsilon,
\end{equation}
and dispersive, in the sense that it satisfies 
bounds similar to the corresponding linear flow as follows, in any time interval $I$ of size
\[
\epsilon^{-4} \lesssim I \ll \epsilon^{-8}
\]
\begin{itemize}
\item global $L^6_{t,x}$ Strichartz bounds, and 
\item global bilinear $L^2_{t,x}$ bounds.
\end{itemize}
\end{conjecture}

Here the exact $L^6_{t,x}$ and $L^2_{t,x}$ bounds are expected
to be optimized exactly by the small solitons 
associated to our problem.

\section{The results so far}
\label{s:results}

There are two settings where the conjectures have been proved until now, namely for semilinear Schr\"odinger flows and for quasilinear Schr\"odinger flows. Other settings are considered in work in progress, particularly for cases where 
the dispersion relation is no longer of Schr\"odinger type.  We discuss these two 
types of problems in what follows.

\subsection{Semilinear Schr\"odinger flows}
The articles \cite{IT-global} and \cite{IT-focusing}
are both devoted to this case, and consider 
the defocusing, respectively the focusing conjectures. This concerns equations of the form
\begin{equation}\label{nls}
i u_t + u_{xx} = C(u,\bar u, u),  \qquad u(0) = u_0, 
\end{equation}
which corresponds to $a(\xi) = \xi^2$.

The symbol $c$ of the trilinear form is assumed to be uniformly smooth, that is,
\begin{equation}\label{c-smooth}
|\partial_\xi^\alpha c(\xi_1,\xi_2,\xi_3)| \leq c_\alpha,
\qquad \xi_1,\xi_2,\xi_3 \in \R,\,   \mbox{  for every  multi-index $\alpha$},
\end{equation}
and conservative,
\begin{equation}\label{c-conserv}
\Im c(\xi,\xi,\eta) = 0, \qquad \xi,\eta \in \R, \mbox { where } \Im z=\mbox{imaginary part of } z\in \mathbb{C}.
\end{equation}
In addition, we consider the defocusing condition
\begin{equation}\label{c-defocus}
c(\xi,\xi,\xi) \geq c > 0, \qquad \xi \in \R  \mbox{ and } c \in \mathbb{R^+}.
\end{equation}
The simplest case here is the classical cubic NLS problem
\begin{equation}\label{nls3}
i u_t + u_{xx} = \pm u|u|^2,  \qquad u(0) = u_0.
\end{equation}
But this particular case has much more structure
than the general case we are considering, and is in particular completely integrable.

\bigskip

The starting point of both works \cite{IT-global} and \cite{IT-focusing} is the following local well-posedness result:

\begin{theorem}[\cite{IT-global}] Under the assumption \eqref{c-smooth}, the equation \eqref{nls}
is locally well-posed in $L^2$.
\end{theorem}

Based on this result, it is natural to 
consider both conjectures for $L^2$ initial data.
The main result of \cite{IT-global} asserts that 
the defocusing conjecture is true in this case.
That was the first global in time well-posedness result of this type.  

\begin{theorem}[\cite{IT-global}]\label{t:main}
Under the above assumptions \eqref{c-smooth}, \eqref{c-conserv} and \eqref{c-defocus}  on the symbol of the cubic form $C$, small initial data 
\[
\|u_0\|_{L^2_x} \leq \epsilon \ll 1,
\] 
yields a unique global solution $u$ for \eqref{nls}, which satisfies the following bounds:

\begin{enumerate}[label=(\roman*)]
\item Uniform $L^2$ bound:
\begin{equation}\label{main-L2}
\| u \|_{L^\infty_t L^2_x} \lesssim \epsilon.
\end{equation}

\item Strichartz bound:
\begin{equation}\label{main-Str}
\| u \|_{L^6_{t,x}} \lesssim \epsilon^\frac23.
\end{equation}

\item Bilinear Strichartz bound:
\begin{equation}\label{main-bi}
\| \partial_x (u \bar u(\cdot+x_0))\|_{L^2_t H_x^{-\frac12}} \lesssim \epsilon^2,
\qquad x_0 \in \R.
\end{equation}
\end{enumerate}
\end{theorem}

Here we note that in the case $x_0 = 0$ the last bound
gives 
\begin{equation}\label{main-bi-diag}
\| \partial_x |u|^2\|_{L^2_tH^{-\frac12}_x} \lesssim \epsilon^2,
\end{equation}
which is the more classical formulation of the bilinear $L^2_{t,x}$ bound. However, 
making this bound uniform  with respect to the  $x_0$ translation
captures the natural separate translation invariance of this  bound, and is also quite useful in our proofs. Even when applied to the classical NLS problem, where $L^2$ global well-posedness was  known before,  this theorem yields a new result:

 \begin{theorem}[\cite{IT-global}]
Consider the defocusing 1-d cubic NLS problem \eqref{nls3}(+) with $L^2$ initial data $u_0$.
Then the global solution $u$ satisfies the following bounds:

\begin{enumerate}[label=(\roman*)]
\item Uniform $L^2$ bound:
\begin{equation}\label{main-L2-model}
\| u \|_{L^\infty_t L^2_x} \lesssim \|u_0\|_{L^2_x}.
\end{equation}

\item Strichartz bound:
\begin{equation}\label{main-Str-model}
\| u \|_{L^6_{t,x}} \lesssim  \|u_0\|_{L^2_x}.
\end{equation}

\item Bilinear Strichartz bound:
\begin{equation}\label{main-bi-model}
\| \partial_x |u|^2\|_{L^2_t (\dot H_x ^{-\frac12} + c L^2_x) } \lesssim  \|u_0\|_{L^2_x}^2, \qquad c = \| u_0\|_{L^2_x}.
\end{equation}
\end{enumerate}
\end{theorem}

Here by scaling one may allow large $L^2_x$ data.
One may compare the above $L^6_{t,x}$ bound with an earlier estimate of Planchon-Vega~\cite{PV}, which applies only to $H^1$ solutions.
\bigskip

The corresponding conjecture for the focusing case was considered in \cite{IT-focusing}:

\begin{theorem}[\cite{IT-focusing}]\label{t:focusing}
Consider the problem \eqref{nls} where the cubic nonlinearity $C$ satisfies \eqref{c-smooth} and \eqref{c-conserv}. Assume that the initial data $u_0$ is small, 
\[
\|u_0\|_{L^2_x} \leq \epsilon \ll 1,
\] 
Then the solution $u$ exists on on a time interval $I_\epsilon = [0,c \epsilon^{-8}]$ and 
has the following properties for every interval $I \subset I_\epsilon$ of size $|I| \leq \epsilon^{-6}$:

\begin{enumerate}[label=(\roman*)]
\item Uniform $L^2$ bound:
\begin{equation}\label{main-L2-focus}
\| u \|_{L^\infty_t(I_\epsilon; L^2_x)} \lesssim \epsilon.
\end{equation}

\item Strichartz bound:
\begin{equation}\label{main-Str-focus}
\| u \|_{L^6_{t,x}(I \times \R)} \lesssim \epsilon^\frac23.
\end{equation}

\item Bilinear Strichartz bound:
\begin{equation}\label{main-bi-focus}
\| \partial_x (u \bar u(\cdot+x_0))\|_{L^2_t (I;H_x^{-\frac12})} \lesssim \epsilon^2,
\qquad x_0 \in \R.
\end{equation}
\end{enumerate}
\end{theorem}
We remark that the intermediate time-scale $\epsilon^{-6}$ does not have an intrinsic meaning from a scaling perspective, but is instead connected to the unit frequency-scale which is implicit in \eqref{c-smooth}.  Under the same hypothesis one could also use a smaller size for $|I|$, in the range $\epsilon^{-4} \lesssim |I| \lesssim \epsilon^{-6}$, with an appropriate adjustment of the constant in \eqref{main-Str-focus}; we refer 
the reader to \cite{IT-focusing} for further details.

\subsection{Quasilinear Schr\"odinger flows}
The aim of our article \cite{IT-qnls} is to prove both conjectures in the setting of quasilinear 
Schr\"odinger flows; this represents the first validation of the conjecture in a quasilinear setting. The are are two types of equations which are considered in  \cite{IT-qnls}. The first is
\begin{equation}\tag{DQNLS}
\label{dqnls}
\left\{ \begin{array}{l}
i u_t + g(u,\partial_x u ) \partial_x^2 u = 
N(u,\partial_x u) , \quad u:
\R \times \R \to \C ,\\ \\
u(0,x) = u_0 (x),
\end{array} \right. 
\end{equation}
with a metric $g$ which is a real valued positive function and a source term $N$ which is a complex valued (real) smooth function of its arguments.
Secondly  we consider the simpler problem
\begin{equation}\tag{QNLS}
\label{qnls}
\left\{ \begin{array}{l}
i u_t + g(u) \partial_x^2 u = 
N(u,\partial_x u) , \quad u:
\R \times \R \to \C, \\ \\
u(0,x) = u_0 (x),
\end{array} \right. 
\end{equation}
where $N$ is  at most quadratic in $\partial u$. This can be seen as the differentiated form of \eqref{dqnls}. For this reason, the results 
for \eqref{qnls} and \eqref{dqnls} are essentially identical, with the only difference that \eqref{dqnls} requires an extra derivative for the 
solutions compared to \eqref{qnls}.

Both of these equations are considered in the cubic setting, which is to say that $g$ is at least quadratic and $N$ is at least cubic. The starting point for the global/long time results is again a corresponding local well-posedness result, which is stated as follows:

\begin{theorem}[\cite{IT-qnls}]\label{t:local}
The cubic problem \eqref{qnls} is locally well-posed for small data in $H^s$ for $s > 1$, and the cubic problem \eqref{dqnls} is locally well-posed for small data in $H^s$ for $s > 2$.
\end{theorem}

The remarkable fact about this result is that it 
uses the tools developed for the global/long time results, and which are discussed in the next section,
in order to drastically improve the the local theory, developed earlier in \cite{KPV} and \cite{MMT2}, by at least one unit and all the way to the sharp threshold. 
Indeed, it is not so difficult to show that under 
the assumptions of the above theorem, \eqref{qnls} 
is generically ill-posed in $H^s$ for $s < 1$, while 
\eqref{qnls} is generically ill-posed in $H^s$ for $s < 2$. 

While the local well-posedness result does not 
require any qualitative assumptions on $g$ and $N$,
for the long time results we assume also that 
\eqref{qnls}/\eqref{dqnls} have the phase rotation symmetry and are conservative in the sense of Definition~\ref{d:conservative}. For the defocusing 
property we use a quantitative bound which reads as 
follows:

\begin{definition}
We say that the equation \eqref{qnls}/\eqref{dqnls}
is \emph{defocusing} if  
    \begin{equation}\label{defocusing-qnls}
    c(\xi,\xi,\xi) \gtrsim \la \xi \ra^{2+2k},
    \end{equation}
    where $k = 1$ for \eqref{dqnls}
    respectively $k=0$ for \eqref{qnls}.
\end{definition}

Now we can state our main results. For the defocusing problem we obtain global solutions, 
validating Conjecture~\ref{c:defocusing}:

\begin{theorem}[\cite{IT-qnls}]\label{t:global}
 a) Assume that the equation \eqref{qnls} has phase rotation symmetry and is conservative and defocusing. Given $s > 1$, assume  that
  the initial data $u_0$ is small in $H^s$,
 \begin{equation}
 \|u_0\|_{H^s_x} \leq \epsilon \ll 1.   
 \end{equation}
Then the solutions are global in time, and satisfy 
\begin{equation}
\| u\|_{L^\infty_t H^s_x} \lesssim \epsilon,
\end{equation}
as well as appropriate $L^6_{t,x}$ Strichartz
and bilinear $L^2_{t,x}$ bounds.

b) Assume that the equation \eqref{dqnls} has phase rotation symmetry and is conservative and defocusing.
 Given $s > 2$, assume  that
the initial data $u_0$ is small in $H^s$,
 \begin{equation}
 \|u_0\|_{H^s_x} \leq \epsilon \ll 1.   
 \end{equation}
Then the solutions are global in time, and satisfy 
\begin{equation}
\| u\|_{L^\infty_t H^s_x} \lesssim \epsilon,
\end{equation}
as well as appropriate $L^6_{t,x}$ Strichartz
and bilinear $L^2_{t,x}$ bounds.
\end{theorem}

For completeness, we describe the  $L^6_{t,x}$ Strichartz and bilinear $L^2_{t,x}$ bounds for \eqref{qnls}; the \eqref{dqnls} bounds are similar but one derivative higher. These are as follows:
\begin{enumerate}[label=(\roman*)]
    \item $L^6_{t,x}$ Strichartz bounds:
  \begin{equation}
   \| \la D\ra^{\frac56+} u\|_{L^6_{t,x}} \lesssim \epsilon .  
  \end{equation}  
\item{Bilinear $L^2_{t,x}$ bounds}, which can be stated in a balanced form
\begin{equation}
   \| \partial |\la D\ra^{\frac34+} u|^2\|_{L^2_{t,x}} \lesssim \epsilon^2   ,
  \end{equation}  
and in an imbalanced form\footnote{ Here $T_u$ is the standard paraproduct operator, see e.g. \cite{Metivier}.}
\begin{equation}
   \| T_u \bu\|_{L^2_tH^{\frac32+}_x} \lesssim \epsilon^2.     \end{equation} 

\end{enumerate}

Finally, we turn our attention to the focusing case:

\begin{theorem}[\cite{IT-qnls}]\label{t:long}
 a) Assume that the equation \eqref{qnls} has phase rotation symmetry and is conservative. Given $s > 1$, assume  that the initial data $u_0$ is small in $H^s$,
 \begin{equation}
 \|u_0\|_{H^s_x} \leq \epsilon \ll 1.   
 \end{equation}
Then the lifespan of the solutions is at least $O(\epsilon^{-8})$.

b) Assume that the equation \eqref{dqnls} has phase rotation symmetry and is conservative. Given $s >2$, 
assume that the initial data $u_0$ is small in $H^s$,
 \begin{equation}
 \|u_0\|_{H^s_x} \leq \epsilon \ll 1.   
 \end{equation}
Then the lifespan of the solutions is at least $O(\epsilon^{-8})$.
\end{theorem}

As in the semilinear case, we also establish the same $L^6_{t,x}$ Strichartz bounds and the bilinear $L^2_{t,x}$ 
bounds on shorter time intervals; the reader is referred to \cite{IT-qnls} for the details.

One final remark is that both theorems above 
hold at exactly the same regularity level as the 
local well-posedness result, which is in turn sharp !

\section{Key ideas and methods}
\label{s:methods}

We begin our discussion with a short list of five  key ideas which play the leading role in our 
approach to these two conjectures: 

\begin{enumerate}
    \item A \emph{bootstrap argument} carried out using \emph{frequency envelopes}, which are used in order 
    to accurately track frequency localized components of the solution over long time scales.
 \medskip
 
 \item \emph{Energy estimates,}  but developed at a local level in the form of  \emph{density-flux identities}.  These are carried out in a nonlocal setting, where both the densities and the fluxes involve translation invariant multilinear forms.

\medskip

\item \emph{Modified energies}, in a manner which  is akin to the \emph{I-method}.
But we  implement this idea in a frequency localized setting and  at the level of density-flux identities, rather than directly for energy functionals.
\medskip

\item \emph{Interaction Morawetz bounds}, also carried out in a frequency localized fashion and
 extended to the setting and language of nonlocal multilinear forms.
\medskip

\item \emph{Strichartz estimates} also play a role  in the study of local well-posedness. These are developed via \emph{wave packet parametrices}, after peeling off ``perturbative'' errors.
\medskip

\end{enumerate}

For the remainder of this section we expand on each of these ideas as well as other related techniques.

\subsection{The Littlewood-Paley decomposition}
This is a standard tool in nonlinear dispersve equations, whose use in this context 
is motivated by a dichotomy arising when considering estimates for multilinear forms:
\begin{itemize}
    \item parallel interactions, corresponding to nearby frequencies with close group velocities.
Estimates for these strongly rely  on $L^6_{t,x}$ Strichartz estimates.
\item transverse interactions, corresponding to separated frequencies and corresponding group velocities. Bounds for such interactions rely instead on bilinear $L^2_{t,x}$ bounds.
\end{itemize}

A typical Littlewood-Paley decomposition in one dimension may be written as 
\[
u = \sum u_k,
\]
where each $u_k$ roughly represents the projection of $u$ to a suitable  Littlewood-Paley region.
The size of these regions has to be carefully tailored to the problem at hand, but the aim of the Littlewood-Paley decomposition remains the same, namely 
(i) to estimate each piece $u_k$ individually, and (ii) to estimate bilinear interactions of separated pieces.

For the semilinear Schr\"odinger problem studied in \cite{IT-global} and \cite{IT-focusing} there are no linear interactions between regions, so we  
are able to use a finer lattice partition of the Fourier space, i.e. into unit size intervals indexed by an integer $k$.

For the quasilinear Schr\"odinger problem studied in \cite{IT-qnls}, on the other hand, the size of the Littlewood-Paley regions is dictated by the maximal frequency spreading allowed by the associated  Hamilton flow.
But this corresponds exactly to the dyadic Littlewood-Paley decomposition, where the summation  index $k$ is a nonnegative integer which corresponds to frequency regions of size $2^k$.
\medskip

\subsection{Frequency envelopes}
This is a tool originally introduced by
Tao~\cite{Tao-WM} in order  to track the time evolution of the energies of each of the Littlewood-Paley pieces in nonlinear evolutions.

In our setting, the strategy is to start with a frequency envelope $\{c_k\} \in \ell^2$ 
which controls the size of the Littlewood-Paley pieces of the initial data, in the sense that
\begin{equation}
\|u_{0k} \|_{H^s_x} \lesssim \epsilon c_k, 
\qquad \|c_k\|_{\ell^2} \lesssim 1.
\end{equation}
Then one would like to show that similar bounds carry over to Littlewood-Paley pieces of the solutions, 
\begin{equation}\label{fe-en}
\|u_{k} \|_{L^\infty_t H^s_x} \lesssim \epsilon c_k.
\end{equation}  
A difficulty one encounters while doing this is that 
there is always leakage between nearby frequencies,
which would defeat this idea if the nearby $c_k$'s 
are completely uncorrelated. This lead to the idea of introducing a key assumption on $c_k$, namely that they should be \emph{slowly varying}. 

In the context of the dyadic Littlewood-Paley decomposition, slowly varying simply means that 
\[
\frac{c_j}{c_k} \leq 2^{\delta|j-k|}, \qquad j,k \in \N .
\]
However, this no longer works in the case of the lattice decomposition, and a replacement is needed. The adapted notion of \emph{slowly varying}, developed in \cite{IT-global}, requires instead  that 
\[
M c_k \lesssim c_k,
\]
where $Mc_k$ denotes the associated maximal function.

Using a frequency envelope approach works particularly well in bootstrap arguments. However,
here it is not sufficient to bootstrap only the bound \eqref{fe-en}, instead one has to also bootstrap appropriate $L^6_{t,x}$ and bilinear $L^2_{t,x}$ bounds, which for instance in the simpler case 
of the semilinear Schr\"odinger problem have the form
\[
\| u_k \|_{L^6_{t,x}} \lesssim  \epsilon c_k 
\]
\[
  \|\partial_x ( u_k  \bar u_j(\cdot+h)\|_{L^2_{t,x}} \lesssim 
\epsilon^2 \la k-j\ra^\frac12 c_k  c_j .
\]  
We note that the idea of bootstrapping simultaneously  Strichartz and bilinear bounds
was explored earlier by the authors in the Benjamin-Ono equation context, see \cite{IT-BO}.

\subsection{A collection of related equations} 
Particularly in the quasilinear context of 
\cite{IT-qnls}, in order to better understand the nonlinear dynamics it is important to separate 
the roles played by different frequency interactions.
For clarity we discuss this in the context of \eqref{qnls}, which we recall here:
\begin{equation}\tag{QNLS}
i u_t + g(u) \partial_x^2 u = 
N(u,\partial_x u). 
\end{equation}
As always, a leading role is played by the linearized equation
\begin{equation}\tag{QNLS-lin}
\label{qnls-lin}
i v_t + g(u) \partial_x^2 v = 
N^{lin}(u) v. 
\end{equation}
But at the heart of both flows above lies the 
associated paradifferential equation, see \cite{Bony,Metivier}, which isolates 
the quasilinear interaction mode, 
\begin{equation}\tag{QNLS-para}
\label{qnls-para}
i w_{k t} + \partial_x g(u_{<k}) \partial_x w_k = 
f_k ,
\end{equation}
and may be thought of as modeling high frequency waves moving on a low frequency background.

One may then rewrite the full equation in paradifferential form, as an infinite system
\begin{equation}\tag{QNLS}
\label{qnls+}
i u_{k t} + \partial_x g(u_{<k}) \partial_x u_k = 
N_k(u,\partial_x u).  
\end{equation}
The same can be done with the linearized equation,
\begin{equation}\tag{QNLS-lin}
\label{qnls=lin+}
i v_{k t} + \partial_x g(u_{<k}) \partial_x v_k = 
N^{lin}_k (u) v. 
\end{equation}
By Bony's formalism, one may expect the source terms on the right to play a perturbative role. This is indeed the case on short time-scales. However, for the long time analysis it is better to separate 
the doubly resonant interactions,
\begin{equation}\tag{QNLS}
\label{qnls++}
i u_{k t} + \partial_x g(u_{<k}) \partial_x u_k = 
N_k^{{{tr}}}(u,\partial_x u)  + C_k(u,\bar u,u). 
\end{equation}
Here $C_k(u,\bar u,u)$ can be thought of as a semilinear term but which carries strong long 
time interactions of balanced frequencies, whereas
the remainder $N_k^{{{tr}}}(u,\partial_x u)$
only carries transversal interactions and is 
fully perturbative.

\subsection{Conservation laws in density-flux form}
Energy estimates play a key role in tracking the 
long time behavior of solutions. In the linear Schr\"odinger context, three key conservation laws 
are provided by the mass, momentum and energy, 
\[
\bM(u) = \int |u|^2 \,dx, \qquad 
\bP(u) = 2 \int \Im ( u \partial_x \bar u) \,dx,
\qquad 
\bE(u) = 4 \int |\partial_x u|^2\, dx. 
\]
These are no longer conserved in our nonlinear flows;
for instance in the case of the mass we obtain a relation of the form
\[
 \frac{d}{dt} {\bf M} = \int C^4_m(u,\bar u, u,\bar u) \, dx
\]
for a quartic (and higher) form $C^4_m$, and similarly for the momentum.

But a more descriptive way to write such conservation  laws is as density-flux identities. 
To achieve this one needs first mass/momentum/energy densities
\[
{\bf M} = \int M(u,\bar u) \, dx, \qquad 
{\bf P} =\int P(u,\bar u) \, dx, \qquad 
{\bf E} =\int E(u,\bar u) \, dx.
\]
These are not uniquely determined, and must be carefully chosen. But if this is done, in the linear case one has exact density-flux identities
\[
 \frac{d}{dt}  M = \partial_x P, \qquad \frac{d}{dt}P = \partial_x E.
\]
In the nonlinear case 
these density-flux identities are no longer exact, 
but instead have the form, say for \eqref{qnls},
\[
\partial_t M(u,\bar u) = \partial_x[g  P(u,\bar u)] + C^4_m(u,\bar u,u,\bar u) + \text{higher},
\]
respectively
\[
\partial_t P(u,\bar u) = \partial_x [ g E(u,\bar u)]  + C^4_p(u,\bar u,u,\bar u)  + \text{higher}.
\]
where ``higher" stands for terms which are at least $6$-linear forms. We use such identities not directly, but instead in a frequency localized form, which is roughly as follows:
\[
\partial_t M_k(u,\bar u) = \partial_x [ g_{<k} P_k(u,\bar u)] + C^4_{m,k}(u,\bar u,u,\bar u) 
+ \text{higher},
\]
\[
\partial_t P_k(u,\bar u) = \partial_x [ g_{<k} E_k(u,\bar u)]  
+ C^4_{p,k}(u,\bar u,u,\bar u) 
+ \text{higher}.
\]
Here for simplicity one could take $M_k(u,\bar u)=
M(u_k,\bar{u}_k)$ and similarly for the momentum.

\subsection{Energy corrections for long time results}
The above density-flux relations are useful for the study of short time dynamics, but not so much for long time dynamics, where quartic forms are not easily estimated directly. However, assuming that 
our equations have the phase rotation symmetry and are conservative, the multilinear forms $ C^4_{m,k}$
and $ C^4_{p,k}$ turn out to have a better structure,
in that their symbols vanish on the doubly resonant set. This property allows us to construct quartic energy corrections of the form
\[
  M_k^\sharp(u,\bar u) = M_k(u,\bar u) + B_{l,m}^4(u,\bar u,u,\bar u),
\]
\[
  P_k^\sharp(u,\bar u) = P_k(u,\bar u) + B_{k,p}^4(u,\bar u,u,\bar u),
\]  
which satisfy better conservation laws
\[
\partial_t M_k^\sharp = \partial_x (g_{<k} P_k +
R_{k,m}^4) + F_{k,m}^{4,tr}   +  R_{k,m}^6 ,
\]
\[
\partial_t P_k^\sharp = \partial_x (g_{<k} E_k +
R_{k,p}^4)  + F_{k,p}^{4,tr}   +  R_{k,p}^6 ,
\]
which is exact modulo quartic terms $F_{k,m}^{4,tr}$ and $F_{k,p}^{4,tr}$ with only transversal interactions, and sixth order terms.

This is somewhat similar to   
the second generation I-method \cite{I-method},
\cite{I-method2}.
We note that constructing these energy corrections requires solving a nontrivial division problem at the symbol level, of the form
\[
c^4 = \Delta^4 \xi^2 \cdot  b^4 + \Delta^4 \xi \cdot  r^4 +
(\xi_{odd}-\xi_{even})^2 q^{4}.
\]

Finally, we remark that the frequency localized uniform energy bounds follow by direct integration
from these identities, as the six-linear error can be estimated from the bootstrap assumptions.

\subsection{Bilinear $L^2$ estimates and Interaction Morawetz identities}
Unlike the energy estimates, bootstrapping the bilinear $L^2_{t,x}$ bounds cannot be done directly 
using linear theory, as (i) the problem is quasilinear and (ii) the nonlinearity is 
nonperturbative.

Instead, for this we rely on   Interaction Morawetz identities, an idea first introduced by the I-team    in the study of the energy critical  NLS problem in
three space dimensions ~\cite{MR2415387}. However, our approach is closer to the one dimensional version of  Planchon-Vega~\cite{PV}. We distinguish two cases:

\medskip

\textbf{A. The diagonal case.} This starts by introducing an  Interaction Morawetz functional for the diagonal case, namely
 \[
I(u_k,u_k) =  \int_{x < y} M_k^\sharp(x) P_k^\sharp(y) -   M_k^\sharp(y)
P_k^\sharp(x) \, dx dy,
 \]
which aims to quantify the self-interaction potential of $u_k$ up to time $t$.
The key point is that its time derivative 
is positive definite at the leading order. Precisely, in the defocusing case we have a relation of the form
 \[
 \frac{d}{dt} I(u_\lambda,u_\lambda) \approx \| \partial_x (u_\lambda \bar
 u_\lambda)\|_{L^2_{t,x}}^2 + \| u_\lambda\|_{L^6_{t,x}}^6 + \text{Errors (6,8,10)},
 \] 
where the errors of varying homogeneities have the common feature that they may be estimated perturbatively using the bootstrap assumptions.

The above relation can be  used to prove the global in time $L^6_{t,x}$ Strichartz bound  and the diagonal bilinear $L^2_{t,x}$ bound, in the defocusing case. In the focusing case the $L^6_{t,x}$ norm 
appears with the wrong sign, and therefore has to 
be moved to the perturbative box and estimated directly using interpolation and H\"older's inequality in time. The last step leads exactly to the $\epsilon^{-8}$ lifespan bound.

\medskip

\textbf{B. The transversal case.}
The transversal  Interaction Morawetz functional
has the form
 \[
I(u_k,u_j) =  \int_{x < y} M_k^\sharp(x) P_j^\sharp(y) -   M_j^\sharp(y)
P_k^\sharp(x) \,dx dy,
 \]
 and its time derivative has the form
\[
 \frac{d}{dt} I(u_\lambda,u_\mu) \approx \| \partial_x (u_\lambda \bar
 u_\mu)\|_{L^2_{t,x}}^2  + \text{Errors (6,8,10)} .
 \] 
 Now all the errors can be estimated perturbatively, 
allows us to prove the off-diagonal bilinear $L^2_{t,x}$ bounds.

\subsection{Strichartz estimates}
While in the case of the semilinear problem 
the frequency localized $L^6_{t,x}$ bounds come from
the diagonal case of the Interaction Morawetz analysis, for the quasilinear case in \cite{IT-qnls} there 
are two distinct $L^6_{t,x}$ bounds which appear in the 
analysis:

\begin{itemize}
    \item A long time $L^6_{t,x}$ bound for the solutions 
    to the nonlinear problem, with a loss of $1/6-$ derivatives.

    \item A short time $L^6_{t,x}$ bound which applies both to solutions for the nonlinear problem
    and the linearized equation, without loss of derivative.
\end{itemize}

The long time $L^6_{t,x}$ bounds come from the Interaction Morawetz identities, and  were discussed in the previous subsection. However, because of the loss of derivatives and lack of a linearized version of such estimates, they cannot be used to close the linear well-posedness.
This is why we need the short time bounds.
These are established at the level of the paradifferential equation
\eqref{qnls-para}.
The main challenge is that we have a variable coefficient problem. The key points of our approach are only briefly discussed below, as they are not directly related to the two conjectures. 

\begin{enumerate}[label=(\roman*)]
\item We flatten the metric with change of coordinates, a strategy previously used by 
Burq-Planchon~\cite{BP} in the case of time independent coefficients. 
    \item We use the equation for $u$ to calculate
    the coefficient of the first order term which represents the contribution arising from the time derivative of the coefficients, and which is overall nonperturbative.
    \item We split the above coefficient  into a high frequency and a low frequency part.
    \item We use bilinear $L^2_{t,x}$ estimates to 
    estimate perturbatively the contribution of the 
    high frequency part of the coefficient. 
    This is achieved by constructing a new, larger space  of ``perturbative" source terms $f_\lambda$ for the paradifferential equation. 
    \item We use the wave packet parametrix of  (Marzuola-Metcalfe-Tataru~\cite{MMT-param}
    to prove the Strichartz estimates for the remaining part of the equation.
\end{enumerate}


\bibliographystyle{hplain}

\end{document}